\documentclass[a4paper,11pt]{article}
\usepackage{graphicx,multirow}
\usepackage{amssymb,amsmath,amsthm}
\usepackage{lscape,lipsum,microtype}
\usepackage[margin=2.54cm]{geometry}  
\usepackage{hyperref}

\usepackage{ifluatex}
\ifluatex
 \usepackage{unicode-math}
\fi

\usepackage{latexsym}
\usepackage{amsmath}
\usepackage{amsfonts}
\usepackage{amssymb}
\usepackage{amscd}

\newcommand{\D }{\Delta }

\newcommand{\e }{\varepsilon }

\newcommand{\n }{\nabla }

\newcommand{\Sig }{\Sigma}

\newcommand{\ov}{\overline}
\newcommand{\intbar}{\mathop{\int\makebox(-13.5,0){\rule[4pt]{.7em}{0.3pt}}%
\kern-6pt}\nolimits}

\newcommand{\be}{\begin{equation}}
\newcommand{\ee}{\end{equation}}
\newcommand{\bes}{\begin{equation*}}
\newcommand{\ees}{\end{equation*}}
\newcommand{\ba}{\begin{eqnarray}}
\newcommand{\ea}{\end{eqnarray}}
\newcommand{\bas}{\begin{eqnarray*}}
\newcommand{\eas}{\end{eqnarray*}}
\newenvironment{pf}{\noindent{\sc Proof}.\enspace}{\rule{2mm}{2mm}\medskip}
\newenvironment{pfn}{\noindent{\sc Proof}}{\rule{2mm}{2mm}\medskip}

\newcommand{\R}{\mathbb{R}}

\newcommand{\N}{\mathbb{N}}

\author{ Cheikh Birahim NDIAYE}

\date{}

\title{\bf Optimal control for the Paneitz obstacle problem}
\begin{document}

\newtheorem{lem}{Lemma}[section]
\newtheorem{pro}[lem]{Proposition}
\newtheorem{thm}[lem]{Theorem}
\newtheorem{rem}[lem]{Remark}
\newtheorem{cor}[lem]{Corollary}
\newtheorem{df}[lem]{Definition}

\maketitle

\begin{center}
{\small

\noindent  Department of Mathematics Howard University \\  Annex $3$, Graduate School of Arts and Sciences, $217$\\ DC 20059 Washington, USA

}

\end{center}

\footnotetext[1]{E-mail addresses: cheikh.ndiaye@howard.edu\\
\thanks{\\ The author was partially supported by NSF grant DMS--2000164.}}

\

\

\begin{center}
{\bf Abstract}
\end{center}
In this paper, we study a natural optimal control problem associated to the Paneitz obstacle problem on closed \;$4$-dimensional Riemannian manifolds. We show the existence of an optimal control which is an optimal state and induces also a conformal metric with prescribed \;$Q$-curvature. We show also \;$C^{\infty}$-regularity of optimal controls and some compactness results for the optimal controls. In the case of the \;$4$-dimensional standard sphere, we characterize all optimal controls. 
 \begin{center}


\bigskip

\centerline{\bf AMS subject classification: 53C21, 35C60, 58J60, 55N10.}

\end{center}


\section{Introduction and statement of the results}
One of the most important problem in conformal geometry is the problem of finding conformal metrics with a prescribed curvature quantity. An example of curvature quantity which has received a lot of attention in the last decades is the Branson's \;$Q$-curvature. It is a Riemannian scalar invariant introduced by Branson-Oersted\cite{bo} (see also Branson\cite{bran1}) for closed four-dimensional Riemannian manifolds. 
\vspace{8pt}

\noindent
Given \;$(M, g)$\; a four-dimensional closed Riemannian manifold with Ricci tensor \;$Ric_{g}$, scalar curvature \;$R_{g}$, and Laplace-Beltrami operator \;$\D_g$, the $Q$-curvature of \;$(M, g)$\; is defined by
\begin{equation}\label{eq:P}
Q_{g}=-\frac{1}{12}(\D_{g}R_{g}-R_{g}^{2}+3|Ric_{g}|^{2}).
\end{equation}
Under the conformal change of metric $g_u=e^{2u}g$\; with \;$u$\; a smooth function on \;$M$, the $Q$-curvature transforms in the following way
\begin{equation}\label{eq:P1}
P_gu+2Q_g=2Q_{g_{u}}e^{4u},
\end{equation}
where \;$P_g$ \;is the Paneitz operator introduced by Paneitz\cite{p1} and is defined by the following formula
\ba\label{eq:P2}
P_{g}\varphi=\D_{g}^{2}\varphi+div_{g}\left((\frac{2}{3}R_{g}g-2Ric_{g})\n_g\varphi\right),
\ea
where \;$\varphi$\; is any smooth function on \;$M$, $div_g$\; is the divergence of with respect to \;$g$, and \;$\n_g$\; denotes the covariant derivative with respect to \;$g$. When, one changes conformally \;$g$\; as before, namely by \;$g_u=e^{2u}g$\; with \;$u$\; a smooth function on \;$M$, $P_g$\; obeys the following simple transformation law
\begin{equation}\label{eq:tranlaw}
P_{g_u}=e^{-4u}P_{g}.
\end{equation}
The equation \eqref{eq:P1} and the formula \eqref{eq:tranlaw} are analogous to classical ones which hold on closed Riemannian surfaces. Indeed, given a closed Riemannian surface $(\Sig, g)$ and \;$g_u=e^{2u}g$\; a conformal change of \;$g$\; with \;$u$\; a smooth function on \;$\Sig$, it is well know that
\ba\label{eq:g}
\Delta_{g_u}=e^{-2u}\Delta_{g},\;\;\;\;\;\;\;\;\;-\Delta_{g}u+K_{g}=K_{g_u}e^{2u},
\ea
where for a background metric \;$\tilde g$\; on \;$\Sigma$ , \;$\Delta_{\tilde g}$\; and \;$K_{\tilde g}$\; are respectively the Laplace-Beltrami operator and the Gauss curvature of \;($\Sigma, \tilde g$). In addition to these, we have an analogy with the classical Gauss-Bonnet formula
$$
\int_{\Sigma}K_{g}dV_{g}=2\pi\chi(\Sigma),
$$
where\;$\chi(\Sigma)$\; is the Euler characteristic of \;$\Sigma$\; and \;$dV_g$\; is the volume form of \;$\Sig$\; with respect to \;$g$. In fact, we have the Chern-Gauss-Bonnet formula
$$
\int_{M}(Q_{g}+\frac{|W_{g}|^{2}}{8})dV_{g}=4\pi^{2}\chi(M),
$$
where \;$W_{g}$\; denotes the Weyl tensor of \;($M,g$) and \;$\chi(M)$\; is the Euler characteristic of  \;$M$. Hence, from the pointwise conformal invariance of  \;$|W_{g}|^{2}dV_{g}$, it follows that \;$\int_{M}Q_{g}dV_g$\; is also conformally invariant and will be denoted by \;$\kappa_g$, namely
\begin{equation}\label{eq:kappa}
\kappa_{g}:=\int_{M}Q_{g}dV_g.
\end{equation}
When \;$(M, g)=(\mathbb{S}^4, g_{\mathbb{S}^4})$\; is the \;$4$-dimensional standard sphere, we have
\begin{equation}\label{kgsphere}
\kappa_g=\kappa_{g_{\mathbb{S}^4}}=8\pi^2.
\end{equation}
\vspace{10pt}

\noindent
Of particular importance in Conformal Geometry is the following  Kazdan-Warner type problem. Given a smooth positive function \;$K$ defined on a closed \;$4$-dimensional Riemannian manifold $(M, g)$, under which conditions on \;$K$\; there exists a  Riemannian metric conformal to \;$g$ \;with $Q$-curvature equal to \;$K$ . Thanks to \eqref{eq:P1}, the problem is equivalent to finding a smooth solution of the fourth-order nonlinear partial differential equation
\ba\label{eq:qequation}
P_{g}u+2Q_{g}=2Ke^{4u}\;\;\;\;\;\;in\;\;M.
\ea
Equation \eqref{eq:qequation} is usually refereed to as the prescribed \;$Q$-curvatre equation and has been studied in the framework of Calculus of Variations, Critical Points Theory, Morse Theory and Dynamical Systems, see \cite{cy}, \cite{dm}, \cite{mal1}, \cite{ms}, \cite{nd5}, \cite{nd6}, \cite{nd7}, and the references therein.  
\vspace{8pt}

\noindent
In this paper, we investigate equation \eqref{eq:qequation} in the context of Optimal Control Theory. Precisely, we study the following optimal control problem for the Paneitz obstacle problem
\begin{equation}\label{stateop}
\text{Finding} \;\;\;u_{min}\in H^2_Q(M)\;\;\text{such that}\;\;\; I(u_{min})=\min _{u\in H^2_Q(M)}I(u),
\end{equation}
where
 $$
 I(u)= \, \left<u,u\right>_{g}  - \kappa_{g} \log\left(\int_{M}Ke^{4T_g(u)} dV_g\right),\;\;\;\;u\in H^{2}_Q(M)
$$
with 
$$
 <u,u>_{g}=\int_M\D_gu\D_gvdV_g+\frac{2}{3}\int_{M}R_g\n_g u\cdot\n_gvdV_g-\int_M2Ric_g(\n_g u, \n_gv)dV_g
$$
 $$
 T_g(u)=\arg{\min_{v\in H^2_Q(M),\\\;\;v\geq u}\left<v, v \right>_g}
 $$
 and 
 $$H^2_Q(M):=\{u\in H^2(M):\;\; \int_MQ_gudV_g=0\}$$
 with \;$H^2(M)$\; denoting the space of functions on \;$M$\; which are of class \;$L^2$, together with their first and second derivatives. Moreover, the symbol $$\arg{\min_{v\in H^2_Q(M),\\\;\;v\geq u}\left<v, v \right>_g}$$ denotes the unique solution to the minimization problem \;$$\min_{v\in H^2_Q(M),\\\;\;v\geq u}\left<v, v \right>_g,$$ see Lemma \ref{obstacleq}. We remark that for \;$u$\; smooth,
 $$
 \left<u,u\right>_{g}=\left<P_gu, u\right>_{L^2(M)},
 $$
where \;$\left<\cdot, \cdot\right>_{L^2(M)}$\;denotes the \;$L^2$\; scalar product.
\vspace{6pt}

\noindent
 In the subcritical case, namely \;$0<\kappa_g<8\pi^2$, we  prove the following result.
\begin{thm}\label{positivevarmass}
Assuming that \;$P_g \geq 0$, $\ker P_g\simeq\R$, and \;$0<\kappa_g <8\pi^2$, then \;there exists \;$$u_{min}\in C^{\infty}(M)\cap H^2_Q(M)$$\; such that
$$
I(u_{min})=\min _{v\in H^2_Q(M)}I(v)\;\;\;\;{and}\;\;\;\;u_{min}= T_g(u_{min}).
$$
Moreover, setting $$
u_c=u_{min}-\frac{1}{4}\log\int_MKe^{4u_{min}}+\frac{1}{4}\log \kappa_g\;\;\;\text{and}\;\;\;\;g_c=e^{2u_c}g,
$$ we have 
$$
Q_{g_c}=K.
$$
\end{thm}
\vspace{10pt}

\noindent
To state our existence result in the critical case, i.e $\kappa_g=8\pi^2$, we first set some notations.
We define \;$\mathcal{F}_K: M\longrightarrow \R$\; as follows
\begin{equation}\label{eq:limitfs}
\mathcal{F}_K(a):=2\left(H(a, a)+\frac{1}{2}\log(K(a))\right),\;\;\;a\in M
\end{equation}
where \;$H$\; is the regular part of the Green's function \;$G$\; of $P_g(\cdot)+2Q_g$ satisfying the  normalization $\int_M Q_g(x) G(\cdot, x)dV_g(x)=0$, see Section \ref{eq:notpre}. Furthermore, we define
\begin{equation}\label{eq:critfk}
Crit(\mathcal{F}_K):=\{a\in M: \;\;a\;\;\text{is critical point of} \;\;\mathcal{F}_K\}.
\end{equation}
Moreover, for \;$a\in M$\; we set
\begin{equation}\label{eq:partiallimit}
\mathcal{F}^a(x):=e^{(H(a, x)++\frac{1}{4}\log(K(x))},\;\;\;x\in M
\end{equation}
and define
\begin{equation}\label{eq:defindexa}
\mathcal{L}_K(a):=- \mathcal{F}^{a}(a)L_g(\mathcal{F}^{a})(a),
\end{equation}
where $$L_g:=-\D_g+\frac{1}{6}R_g$$ is the conformal Laplacian associated to \;$g$.  We set also
\begin{equation}\label{eq:critsett}
\mathcal{F}_{\infty}^+:=\{a\in Crit(\mathcal{F}_K):\:\;\mathcal{L}_K(a)>0\}.
\end{equation}
With this notation, our existence result in the critical case reads as follows:
\begin{thm}\label{positivevarmassi}
Assuming that \;$P_g \geq 0$, $\ker P_g\simeq\R$,  \;$\kappa_g =8\pi^2$, and \;$\mathcal{F}_{\infty}^+=Crit(\mathcal{F}_K)$, then \;there exists \;$$u_{min}\in C^{\infty}(M)\cap H^2_Q(M)$$\; such that
$$
I(u_{min})=\min _{v\in H^2_Q(M)}I(v)\;\;\;\;{and}\;\;\;\;u_{min}= T_g(u_{min}).
$$
Moreover, setting $$
u_c=u_{min}-\frac{1}{4}\log\int_MKe^{4u_{min}}+\frac{1}{4}\log \kappa_g\;\;\text{and}\;\;\;\;g_{c}=e^{2u_c}g ,
$$ 
we have 
$$
Q_{g_c}=K.
$$
\end{thm}
\begin{rem}
\begin{itemize}
\noindent

\item The relation \;$u_{min}=T_g(u_{min})$\; in the above theorems is an additional information with respect to the existence results based on Calculus of Variations, Critical Points Theory, Morse Theory, and Dynamical Systems. It provides the inequality
\begin{equation}\label{obsinq}
\left<u_{min}, u_{min}\right>_g\leq \left<u, u\right>_g, \;\;\;\;\;\forall \;u_{min}\leq  u\in H^2_Q(M).
\end{equation}
\item
We remark that the nonlocal character of \;$e^{4T_g(u)}$\; in the definition of  \;$I$ with respect to  \;$e^{4u}$\; appearing in the definition of  \;$J$\; defined by
$$
J(u):= \, <u,u>_{g}  - \kappa_{g} \log\left(\int_{M}Ke^{u} dV_g\right),\;\;\;\;u\in H^{2}_Q(M)
$$
used in the existence approaches of \eqref{eq:qequation} via Calculus of Variations, Critical Points Theory, Morse Theory, and Dynamical Systems. The trade of the local character to non-local is in contrast with the traditional approach in the study of Differential Equations, but have the advantage of providing automatically the variational inequality \eqref{obsinq}.
\item
The \;$Q$-curvature functional \;$J$\; is invariant by translation by constants, while the \;$Q$-optimal control functional \;$I$\; is not. The functional \;$J$ is weakly lower semicontinuous, but the functional \;$I$\; is not. This makes it difficult to apply the Direct Methods in Calculus of Variations to study \eqref{stateop}.
\item 
We expect formula \eqref{obsinq} to be useful to deal with the case \;$\kappa_g=8\pi^2$\; by helping to track down the loss of coercivity in the Variational Analysis of equation \eqref{eq:qequation}.
\end{itemize}
\end{rem}
\vspace{10pt}

\noindent
As a byproduct of our existence argument, we have the following regularity result for solutions of the optimal control problem \eqref{stateop}.
\begin{thm}\label{regularity}
Assuming that \;$P_g \geq 0$, $\ker P_g\simeq \R$,  \;$0<\kappa_g \leq 8\pi^2$, and \;$u\in H^2_Q(M)$\; is a minimizer of \;$I$\; on \;$H^2_Q(M)$, then
$$
u\in  C^{\infty}(M).
$$
\end{thm} 
\vspace{10pt}

\noindent
An other consequence of our existence argument is the following compactness theorems for the set of minimizers of \;$I$\; on \;$H^2_Q(M)$.  We start with the subcritical case.
\begin{thm}\label{compsub}
Assuming that \;$P_g \geq 0$, $\ker P_g\simeq \R$, and \;$0<\kappa_g <8\pi^2$, then \;$\forall m\in \N$\; there exists \;$C_m>0$\; such that \;$\forall u \in C^{\infty}(M)\cap H^2_Q(M)$\; minimizer of \;$I$\;  on \;$H^2_Q(M)$, we have
$$
||u||_{C^k(M)}\leq C_m.
$$
\end{thm}
\vspace{8pt}

\noindent
For the critical case, setting
\begin{equation}\label{eq:critsett0}
\mathcal{F}_{\infty}^0:=\{a\in Crit(\mathcal{F}_K):\:\;\mathcal{L}_K(a)\neq 0\},
\end{equation}
we have:
\begin{thm}\label{compcri}
Assuming that \;$P_g \geq 0$, $\ker P_g\simeq\R$,  \;$\kappa_g =8\pi^2$, and \;$\mathcal{F}_{\infty}^0=Crit(\mathcal{F}_K)$, then \;$\forall m\in \N$\; there exists \;$C_m>0$\; such that \;$\forall u \in C^{\infty}(M)\cap H^2_Q(M)$\; minimizer of \;$I$\;  on \;$H^2_Q(M)$,  we have
$$
||u||_{C^m(M)}\leq C_m.
$$
\end{thm}
\vspace{10pt}

\noindent
We prove also some results in the particular case of the \;$4$-dimensional standard sphere, see Theorem \ref{sharpomt} and Corollary \ref{sfixedt} in Section \ref {casesphere}.
\vspace{10pt}

\noindent
To prove Theorem \ref{positivevarmass}-Theorem \ref{compcri}, we first use the variational characterization of the solution of  Paneitz obstacle problem \;$T_g(u)$ (see Lemma \ref{obstacleq})\;  to show that the Paneitz obstacle solution map \;$T_g$\; is  idempotent, i.e \;$T^2_g=T_g$, see Proposition \ref{nilpotent}. Next, using the idempotent property of \;$T_g$, we establish some monotonicity formulas, see Lemma \ref{decreasingformula}, Lemma \ref{minimaltunnel}, and Lemma \ref{decreasingformulaop}. Using the later monotonicity formulas, we show that  any minimizer of \;$J$\; or any solution of the optimal control problem \eqref{stateop} is a fixed point of \;$T_g$, see Corollary \ref{rigidityminimizer} and Corollary \ref{rigidityminimizeri}. This allows us to show that the \;$Q$-curvature functional \;$J$\; and the \;$Q$-optimal functional have the same minimizers on \;$H^2_Q(M)$, see Proposition \ref{sameminimizer} . With this at hand, Theorem \ref{positivevarmass} follows from the work of Chang-Yang\cite{cy} in the subcritical case, while Theorem \ref{positivevarmassi} follows from our work in the critical case in \cite{nd6}. Moreover, Theorem \ref{regularity} follows from the regularity result of Uhlenbeck-Viaclosky\cite{uv}. Furthermore, Theorem \ref{compsub} follows from the compactness result of Malchiodi\cite{mal} and Druet-Robert-\cite{dr}, while Theorem \ref{compcri} follows our compactness theorem in \cite{nd6}.
\vspace{10pt}

\noindent
The structure of the paper is as follows. In Section \ref{eq:notpre}, we collect some preliminaries and fix some notations. In Section \ref{opp}, we discuss the Paneitz obstacle problem and some monotonicity formulas involving the \;$Q$-curvature functional \;$J$. We also present some consequences of the latter monotonicity formulas. In Section \ref{socp}, we establish some monotonicity formulas for the \;$Q$-optimal functional \;$I$\; and their consequences as well. In Section \ref{pots}, we present the proof of Theorem \ref{positivevarmass}-Theorem \ref{compcri}. Finally, in Section \ref{casesphere}, we discuss the particular case of the \;$4$-dimensional standard sphere.
\section{Notations and Preliminaries}\label{eq:notpre}
 In this brief section, we fix our notations and give some preliminaries. First of all, from now until the end of the paper, \;$(M, g)$\;  and $\;K: M\longrightarrow \R_+$\; are respectively the given underlying closed four-dimensional Riemannian manifold and the smooth positive function to prescribe.
 \vspace{6pt}
 
  \noindent
We recall the function \;$J$\; used in other approaches to study \eqref{eq:qequation}. 
\begin{equation}\label{eq:defj1}
J(u):= \, \left<u,u\right>_g \, + \, 4 \int_M Q_g u dV_g \,  - \,\kappa_{g} \log\left(\int_{M}Ke^{4u} dV_g\right),\;\;\;\;u\in H^2(M).
\end{equation}
Moreover, we recall the perturbed functional \;$J_t$ \;($0<t\leq 1$) which plays also an important role in the study of  minimizers of \;$J$.
\begin{equation}\label{eq:defjt1}
J_t(u):= \, \left<u,u\right>_g \, + \, 4t \int_M Q_g u dV_g \,  - \,t\kappa_{g} \log\left(\int_{M}Ke^{4u} dV_g\right),\;\;\;\;u\in H^2(M).
\end{equation}
We observe that
$$
J=J_1.
$$
Moreover, we define
$$
\ov{(u)}_Q=\frac{1}{\kappa_g}\int_MQ_gudV_g, \;\;\;\;u\in H^2(M),
$$
so that
$$
H^2_Q(M)=\{u\in H^2(M):\;\;\;\ov{(u)}_Q=0\}.
$$
For \;$a\in M$, we let \;$G (a, \cdot)$\; be the unique solution of the following system
\begin{equation}\label{eq:defG4}
\begin{cases}
P_g G(a, \cdot)+2Q_g(\cdot)=16\pi^2\delta_a(\cdot)\;\;\text{in}\;\;\;\;M\\
\int_M Q_g(x) G(a, x)dV_g(x)=0.
\end{cases}
\end{equation}
It is a well know fact that \;$G(\cdot, \cdot)$\; has a logarithmic singularity. In fact \;$G(\cdot, \cdot)$\; decomposes as follows
\begin{equation}\label{eq:decompG4}
 G(a,x)=\log \left(\frac{1}{\chi^2(d_{g}(a, x))}\right)+H(a, x),  \;\;\;x\neq a \in M.
\end{equation}
where \;$H(\cdot, \cdot)$ is the regular par of \;$G(\cdot, \cdot)$ and  \;$\chi$\; is some smooth cut-off function, see for example \cite{zw}. 
\vspace{8pt}

\noindent
The decomposition of the Green's function \;$G$\; and the arguments of the proof of the Moser-Trudinger's inequality of  Chang-Yang\cite{cy} imply the following Moser-Trudinger type inequality.
\begin{pro}\label{cmt}
Assuming that  \;$P_g \geq 0$, $\ker P_g=\R$,  then \;there exists \; $C=C(M, g)>0$\; such that 
$$
\log \int_Me^{4u}dV_g\leq C+\frac{1}{8\pi^2}\left<u, u\right>_g, \;\;\;\forall u\in H^2_Q(M).
$$
\end{pro}
\vspace{6pt}

\noindent
When \;$(M, g)=(\mathbb{S}^4, g_{\mathbb{S}^4})$, we say \;$v$\; is a standard bubble if 
 
\begin{equation}\label{standard}
P_{g_{\mathbb{S}^4}} v+6=6e^{4v} \;\;\text{on}\;\;\;\mathbb{S}^4.
\end{equation}
By the result of Chang-Yang \cite{cy1}, \;$v$\; satisfies
$$
e^{2v}g_{\mathbb{S}^4}=\varphi^*(g_{\mathbb{S}^4}),.
$$
for some \;$\varphi$\; conformal transformation of \;$\mathbb{S}^4$. It is well-known that the standard bubbles are related to the classical Moser-Trudinger-Onofri inequality. Indeed, we have:
\begin{pro}\label{cmto}
Assuming that \;$(M, g)=(\mathbb{S}^4, g_{\mathbb{S}^4})$\; and \;$K=1$, then
\begin{equation}\label{icmto}
J(u)\geq 0,\;\;\;\;\forall u\in H^2(M).
\end{equation}
Moreover, equality in \eqref{icmto} holds if and only if \;$$v:=u-\frac{1}{4}\log \int _Me^{4u}+\frac{1}{4}\log \frac{\kappa_g}{3}$$ is a standard bubble.
\end{pro}
\vspace{8pt}

\noindent
To end this section, we say \;$w$\; is a \;$Q$-normalized standard bubble, if 
\begin{equation}\label{qnormalized}
w=v-\ov{(v)}_Q,
\end{equation}
with \;$v$\;a standard bubble.
\section{Obstacle problem for the Paneitz operator}\label{opp}
In this section, we study the obstacle problem for the Paneitz operator. Indeed in analogy to the classical obstacle problem for the Laplacian, given \;$u\in H^2_Q(M)$, we look for a solution to the minimization problem 
\begin{equation}\label{obspan}
\min_{v\in H^2_Q(M),\\\;\;v\geq u}\left<v, v \right>_g.\end{equation} We start with the following lemma providing the existence and unicity of solution for the obstacle problem for the Paneitz operator \eqref{obspan}. 
\begin{lem}\label{obstacleq}
Assuming that \;$P_g \geq 0$\;  and \;$\ker P_g \simeq \R$, then \;$\forall u\in H^2_Q(M)$,
there exists a unique $T_g(u)\in H^2_Q(M)$ such that
\begin{equation}\label{tg}
\left<T_g(u),T_g(u) \right>_g=\min_{v\in H^2_Q(M),\\\;\;v\geq u}\left<v, v \right>_g
\end{equation}
\end{lem}
\begin{pf}
Since \;$P_g$\; is self-adjoint, $P_g\geq 0$\; and \;$\ker P_g\simeq\R$, then \;$<\cdot, \cdot>_g$\; defines a scalar product on \;$H^2_Q(M)$\; inducing a norm equivalent to the standard \;$H^2(M)$-norm on \;$H^2_Q(M)$. Hence, as in the classical obstacle problem for the Laplacian, the lemma follows from Direct Methods in the Calculus of Variations.
\end{pf}
\vspace{6pt}

\noindent
We study now some properties of the obstacle solution map \;$T_g\; : H^2_Q(M)\longrightarrow  H^2_Q(M)$. We start with the following algebraic one.
\begin{pro}\label{nilpotent}
Assuming that \;$P_g \geq 0$, \;$\ker P_g \simeq \R$, then the obstacle solution map \;$T_g\; : H^2_Q(M)\longrightarrow  H^2_Q(M)$\; is idempotent, i.e $$T^2_g=T_g.$$
\end{pro}

\begin{pf}
Let $v\in H^2_Q(M)$ such that $v\geq T_g(u)$. Then $T_g(u)\geq u$ implies
$v\geq u$. Thus by minimality, we obtain \;$$\left<v, v \right>_g \geq \left<T_g(u),T_g(u)\right>_g.$$ Hence, since \;$T_g(u)\geq T_g(u)$\; then by unicity we have
$$
T_g(T_g(u))=T_g(u),
$$
thereby ending the proof.
\end{pf}


\vspace{8pt}

\noindent
Next, we discuss some monotonicity formulas. We start with the following one.
\begin{lem}\label{decreasingformula}
Assuming that \;$P_g \geq 0$, \;$\ker P_g \simeq \R$, \;$0<t\leq 1$\; and \;$0<\kappa_g\leq 8\pi^2$, then $$
J_t(u)-J_t(T_g(u))\geq \left<u, u\right>_g-\left<T_g(u)),T_g(u) \right>_g\geq 0, \;\;\;\forall u\in H^2_Q(M).
$$
\end{lem}
\begin{pf}
Using the definition of \;$J_t$ (see \eqref{eq:defjt1}), we have
\begin{equation}
J_t(u)-J_t(T_g(u))= \left<u, u\right>_g-\left<T_g(u),T_g(u) \right>_g-t\kappa_g\left(\log\frac{ \int_M Ke^{4u}dV_g}{ \int_M Ke^{4T_g(u)}dV_g}\right).
\end{equation}
Hence the result follows from \;$K>0$,  \;$T_g(u)\geq u$, and Lemma \ref{obstacleq} .
\end{pf}
\vspace{6pt}

\noindent
Lemma \ref{decreasingformula} imply the following rigidity result.
\begin{cor}\label{rigidity}
Assuming that \;$P_g \geq 0$, $\ker P_g\simeq\R$, \;$0<t\leq 1$\; and \;$0<\kappa_g\leq 8\pi^2$, then \;\;$\forall u\in H^2_Q(M)$,
\begin{equation}\label{ineq}
J_t(T_g(u)) \leq  J_t(u)
\end{equation}
and  
\begin{equation}\label{eql}
J_t(u)=J_t(T_g(u))\implies u = T_g(u).
\end{equation}
\end{cor}
\begin{pf}
Using lemma \ref{decreasingformula}, we have
\begin{equation}\label{eq1}
J_t(u)-J_t(T_g(u))\geq \left<u, u\right>_g-\left<T_g(u),T_g(u) \right>_g\geq 0.
\end{equation}
Thus, \eqref{ineq} follows from \eqref{eq1}. If \;$J_t(u)=J_t(T_g(u))$, then \eqref{eq1} implies
$$
\left<u, u\right>_g=\left<T_g(u),T_g(u) \right>_g.
$$
Hence, since \;$u\geq u$, then the unicity part in Lemma \ref{obstacleq} implies
$$
u=T_g(u),
$$
thereby ending the proof of the corollary.
\end{pf}


\vspace{8pt}

\noindent
Corollary \ref{rigidity} implies that minimizers of \;$J_t$\; on \;$H^2_Q(M)$\; are fixed points of  the obstacle solution map \;$T_g$. Indeed, we have:
\begin{cor}\label{rigidityminimizer}
Assuming that \;$P_g \geq 0$, $\ker P_g\simeq\R$, \;$0<t\leq 1$\; and \;$0<\kappa_g\leq 8\pi^2$, then $$u\in H^2_Q(M)\;\;\;\text{is a minimizer of}\;\;\; J_t\implies u = T_g(u).$$
\end{cor}
\begin{pf}
$u\in H^2_Q(M)$\; is a minimizer of  \;$J_t$\; on \;$H^2_Q(M)$\; implies
\begin{equation}\label{eq3}
J_t(u)\leq J_t(T_g(u)).
\end{equation}
Thus combining \eqref{ineq} and \eqref{eq3}, we get
\begin{equation}\label{eq4}
J_t(u)= J_t(T_g(u)).
\end{equation}
Hence, combining \eqref{eql} and \eqref{eq4}, we obtain
$$
u=T_g(u).
$$
\end{pf}

\begin{rem}
Under the assumption of Corollary \ref{rigidity}, we have Proposition \ref{nilpotent} and Corollary \ref{rigidity} imply that we can assume without loss of generality that any minimizing sequence \;$(u_l)_{l\geq 1}$\; of \;$J_t$\; on  \;$H^2_Q(M)$\;satisfies
$$
u_l=T_g(u_l), \;\;\;\;\forall l\ge 1.
$$
\end{rem}

\section{Optimal control for the Paneitz operator}\label{socp}
In this section, we study a natural optimal control problem associated to the obstacle problem for the Paneitz operator . Indeed, we look for solutions of 
$$\min_{u\in H^2_Q(M)}I(u),$$ where \;$I$\; is the $Q$-optimal control functional defined by

\begin{equation}\label{eq:defi12}
I(u):= \, \left<u,u\right>_g  - \kappa_{g} \log\left(\int_{M}Ke^{4T_g(u)} dV_g\right),\;\;\;\;u\in H^{2}_Q(M).
\end{equation}
Similarly to the \;$Q$-curvature functional $J$, for $0<t\leq 1$\; we define \;$I_t$\; by
\begin{equation}\label{eq:defi12t}
I_t(u):= \, \left<u,u\right>_g  - t\kappa_{g} \log\left(\int_{M}Ke^{4T_g(u)} dV_g\right),\;\;\;\;u\in H^{2}_Q(M).
\end{equation}
\vspace{6pt}

\noindent
We start with the following comparison result. 
\begin{lem}\label{minimaltunnel}
Assuming that \;$P_g \geq 0$, $\ker P_g\simeq\R$, \;$0<t\leq 1$\; and \;$0<\kappa_g\leq 8\pi^2$, then 
$$
I_t\leq J_t\;\;\;\;\;\;\text{on}\;\;\;\;\; H^2_Q(M)\;\;\;\;\;
\text{and}\;\;\;\;\;\;
J_t\circ T_g=I_t\circ T_g\;\;\;\;\text{on}\;\;\; H^2_Q(M).
$$
\end{lem}
\begin{pf}
By definition of \;$J_t$\; and \;$I_t$ (see \eqref{eq:defjt1} and \eqref{eq:defi12t}), we have
$$
J_t(u)-I_t(u)=t\kappa_g\log \left(\frac{\int_{M}Ke^{4T_g(u)}}{\int_{M}Ke^{4u} } \right).
$$
Thus \;$I_t(u)\leq J_t(u)$\; follows from \;$T_g(u)\geq u$\; and \;$K>0$. Moreover, we have
$$
J_t(T_g(u))-I_t(T_g(u))=t\kappa_g\log \left(\frac{\int_{M}Ke^{4T^2_g(u)}}{\int_{M}Ke^{4T_g(u)} }\right) .
$$
Hence, \;$T_g^2=T_g$\; (see Lemma \ref{nilpotent}) implies \;$$J_t(T_g(u))=I_t(T_g(u)).$$
\end{pf}
\vspace{6pt}

\noindent
We have the following monotonicity formula for the \;$Q$-optimal control functional \;$I_t$.
\begin{lem}\label{decreasingformulaop}
Assuming that \;$P_g \geq 0$, \;$\ker P_g \simeq \R$, \;$0<t\leq 1$\; and \;$0<\kappa_g\leq 8\pi^2$, then \;$\forall u\in H^2_Q(M)$,
$$
I_t(u)-I_t(T_g(u))=\left<u, u\right>_g-\left<T_g(u),T_g(u)\right>_g\geq 0.
$$
\end{lem}
\begin{pf}
By definition of \;$I_t$ (see \eqref{eq:defi12t}), we have
$$
I_t(u)-I_t(T_g(u))=\left<u, u\right>_g-\left<T_g(u),T_g(u) \right>_g-t\kappa_g\log \left(\frac{\int_{M}Ke^{4T_g(u)}}{\int_{M}Ke^{4T_g^2(u)} } \right).
$$
Using \;$T_g^2(u)=T_g(u)$\; and the definition of \;$T_g$\; (see Lemma \ref{obstacleq}), we  get
$$I_t(u)-I_t(T_g(u))=\left<u, u\right>_g-\left<T_g(u),T_g(u) \right>_g\geq 0.$$
\end{pf}
\vspace{6pt}

\noindent
Lemma \ref{obstacleq} and Lemma \ref{decreasingformulaop} imply that minimizers of \;$I_t$\; are fixed points of \;$T_g$.
\begin{cor}\label{rigidityminimizeri}
Assuming that \;$P_g \geq 0$, \;$\ker P_g=\R$, \;$0<t\leq 1$\; and \;$0<\kappa_g\leq 8\pi^2$, then  $$u\in H^2_Q(M)\;\;\;\;\;\text{is a minimizer of }\;\;\;\;\; \;I_t\implies u = T_g(u).$$
\end{cor}
\begin{pf}
$u\in H^2_Q(M)$\;is a minimizer of \;$I_t$\; implies 
$$I_t(u)\leq I_t(T_g(u)).$$
Thus Lemma \ref{decreasingformulaop} gives
$$
\left<u, u\right>_g=\left<T_g(u),T_g(u) \right>_g.
$$
Hence, by unicity we have
$$u=T_g(u).$$
\end{pf}
\vspace{6pt}

\noindent
\begin{rem}
Under the assumptions of Corollary \ref{decreasingformulaop}, we have that Proposition \ref{nilpotent} and Corollary \ref{decreasingformulaop} imply that for a minimizing sequence \;$(u_l)_{l\geq 1} $\; of \;$I_t$\; on  \;$H^2_Q(M)$, we can  assume without loss of generality that
$$
u_l=T_g(u_l), \;\;\;\;\forall l\ge 1.
$$
\end{rem}
\vspace{10pt}

\noindent
We have the following proposition showing that \;$I_t$\; and \;$J_t$\; have the same minimizers on \;$H^2_Q(M)$. 
\begin{pro}\label{sameminimizer}
Assuming that \;$P_g \geq 0$, \;$\ker P_g\simeq\R$, \;$0<t\leq 1$\; and \;$0<\kappa_g\leq 8\pi^2$, then $$u\in H^2_Q(M)\;\;\text{is a minimizer of}\;\;J_t\;\;\text{is equivalent to}\;\;u\in H^2_Q(M)\;\;\text{is a minimizer of }\;\;I_t.$$
\end{pro}
\begin{pf}
Suppose $u\in H^2_Q(M)$\; is a minimizer of \;$J_t$. Then Corollary \ref{rigidityminimizer} implies
$$
u=T_g(u).
$$
Thus using Lemma \ref{minimaltunnel} we have
$$
I_t(u)=J_t(u)
$$
For $v\in H^2_Q(M)$, we have Lemma \ref{minimaltunnel}, Lemma \ref{decreasingformulaop}, and $u\in H^2_Q(M)$\; is a minimizer of \;$J_t$\; imply
$$
I_t(v)\geq I_t(T_g(v))=J_t(T_g(v))\geq J_t(u)=I_t(u).
$$
Hence $u\in H^2_Q(M)$\; is a minimizer of \;$I_t$\;on \;$H^2_Q(M)$. Similarly, suppose  $u\in H^2_Q(M)$\; is a minimizer of \;$I_t$. Then Corollary \ref{rigidityminimizeri} implies
$$
u=T_g(u).
$$
Thus using again Lemma \ref{minimaltunnel} we have
$$
I_t(u)=J_t(u).
$$
For $v\in H^2_Q(M)$, we have Lemma \ref{minimaltunnel}, Lemma \ref{decreasingformula}, and \;$u\in H^2_Q(M)$\; is a minimizer of \;$I_t$\; imply
$$
J_t(v)\geq J_t(T_g(v))=I_t(T_g(v))\geq I_t(u)=J_t(u).
$$
Hence \;$u\in H^2_Q(M)$\; is a minimizer of \;$J_t$\; on \;$H^2_Q(M)$.
\end{pf}
\vspace{8pt}

\noindent


\section{Proof of Theorem \ref{positivevarmass} -Theorem \ref{compcri}}\label{pots}
In this section, we present the proof of Theorem \ref{positivevarmass} -Theorem \ref{compcri}. As already mentioned in the introduction, the proofs are based on Proposition \ref{sameminimizer} and some contributions of Chang-Yang\cite{cy}, Druet-Robert\cite{dr}, Malchiodi\cite{mal}, the author\cite{nd6} and Uhlenbeck-Viaclovsky\cite{uv} in the the study of the fourth-order nonlinear partial differential equation \eqref{eq:qequation}.
\vspace{8pt}

\noindent
\begin{pfn}{ of Theorem \ref{positivevarmass}}\\
Since \;$P_g \geq 0$, $\ker P_g=\R$, and \;$0<\kappa_g <8\pi^2$, then the works of Chang-Yang\cite{cy} and Uhlenbeck-Viaclosvky\cite{uv} imply the existence of \;$u_0\in C^{\infty}(M)$ such that
$$
J(u_0)=\min_{u\in H^2(M)} J(u).
$$
Since \;$J$\; is translation invariant, then setting $$u_{min}=u_0-\ov{(u_0)}_Q,$$ we have $$u_{min}\in C^{\infty}(M)\cap H^2_Q(M)$$ and 
$$
J(u_{min})=\min_{u\in H^2_Q(M)} J(u).
$$
Using Proposition \ref{sameminimizer}, we get
$$
I(u_{min})=\min_{u\in H^2_Q(M)} I(u)
$$
Thus Corollary \ref{rigidityminimizeri} implies
$$
u_{min}=T_g(u_{min}).
$$
Recalling that $$J(u_{min})=J(u_0)=\min_{u\in H^2(M)} J(u),$$ we have 
$$
P_gu_{min}+2Q_g=2\kappa_g\frac{Ke^{4u_{min}}}{\int_MKe^{4u_{min}}}.
$$
Thus, setting 
$$
u_c=u_{min}-\frac{1}{4}\log\int_MKe^{4u_{min}}+\frac{1}{4}\log \kappa_g,
$$ we have 
$$
P_gu_{c}+2Q_g=2Ke^{4u_{c}}.
$$
Hence, setting
$$
g_{u_c}=e^{2u_c}g,
$$  
we obtain
$$
Q_{g_{u_c}}=K.
$$
thereby ending the proof.
\end{pfn}

\begin{pfn} {of Theorem \ref{positivevarmassi}}\\
Let \;$\varepsilon_l\in (0, 1)$\; with \;$\varepsilon_l\rightarrow 0$. For \;$l\geq 1$, we define
$$
J_l:=J_{1-\varepsilon_l}\;\;\;\text{and}\;\;\;I_l:=I_{1-\varepsilon_l}
$$
As in the proof of Theorem \ref{positivevarmass}, for \;$l\geq 1$\; the works of Chang-Yang\cite{cy} and Uhlenbeck-Viaclosvky\cite{uv} give the existence of $$u^{l}_{min}\in C^{\infty}(M)\cap H^2_Q(M)$$ such that
\begin{equation}\label{eq:minjl}
J_{l}(u^{l}_{min})=\min_{u\in H^2(M)} J_{l}(u).
\end{equation}
Thus, using Proposition \ref{sameminimizer}, we get
\begin{equation}\label{minoptep}
I_{l}(u^{l}_{min})=\min_{u\in H^2_Q(M)} I_{l}(u).
\end{equation}
Clearly \eqref{eq:minjl} imply, 
\begin{equation}\label{presqlimit}
P_g u^{l}_{min}+2Q_g(1-\varepsilon_l)=2\kappa_g(1-\varepsilon_l)\frac{Ke^{4u^{l}_{min}}}{\int_MKe^{4u^{l}_{min}}}.
\end{equation}
Hence, setting 
\begin{equation}\label{ul}
u_c^{l}=u_{min}^{l}-\frac{1}{4}\log\int_MKe^{4u_{min}^{l}}+\frac{1}{4}\log \kappa_g,
\end{equation}
we obtain
\begin{equation}\label{eulerl1}
P_gu_{c}^{l}+2Q_g(1-\varepsilon_l)=2K(1-\varepsilon_l)\e^{4u_{c}^{l}}.
\end{equation}
Thus our bubbling rate formula  in \cite{nd6} and the assumption \;$\mathcal{F}_{\infty}^+=Crit(\mathcal{F}_K)$\;  prevents the sequence \;$u_{c}^{l}$\; from bubbling. Hence we have
\begin{equation}\label{smoothc}
u_{c}^{l}\longrightarrow u_c \;\;\;\;\text{smoothly, as } \;\;l\longrightarrow \infty.
\end{equation}
Thus \eqref{eulerl1} gives
 \begin{equation}\label{euler1}
 P_gu_{c}+2Q_g=2Ke^{4u_{c}}.
 \end{equation}
Recalling \;$u_{min}^l\in H^2_Q(M)$, we have  \eqref{ul} and \eqref{smoothc} imply
\begin{equation}\label{smoothconv}
u_{min}^{l}\longrightarrow u_{min} \;\;\;\;\text{smoothly}.
\end{equation}
and 
$$
u_c=u_{min}-\frac{1}{4}\log\int_MKe^{4u_{min}}+\frac{1}{4}\log \kappa_g.
$$
Clearly \eqref{smoothconv} and \eqref{minoptep} imply
$$
I(u_{min})=\min_{u\in H^2_Q(M)} I(u).
$$
Hence Corollary \ref{rigidityminimizeri} and \eqref{euler1} imply
$$
u_{min}=T_g(u_{min}).
$$
and 
$$
Q_{g_{u_c}}=K.
$$
\end{pfn}
\vspace{6pt}

\noindent
\begin{pfn}{ of Theorem \ref{regularity}}\\
It follows  directly from Proposition \ref{sameminimizer}, the translation invariant property of \;$J$\; and the regularity result of Uhlenbeck-Viaclovsky\cite{uv}.
\end{pfn}
\vspace{6pt}

\noindent
\begin{pfn}{ of Theorem \ref{compsub}}\\
Let $u\in C^{\infty}(M)\cap H^2_Q(M)$ be a minimizer of \;$I$\; on \;$H^2_Q(M)$. Then  the translation invariance property of \;$J$\; and Proposition \ref{sameminimizer} imply \;$u$\; is a minimizer of \;$J$\; on \;$H^2(M)$. Hence \;$u$\; satisfies
$$
P_gu+2Q_g=2\kappa_g \frac{Ke^{4u}}{\int_MKe^{4u}}.
$$
Then, setting 
\begin{equation}\label{defiv}
v=u-\frac{1}{4}\log \int_{M}Ke^{4u}+\frac{1}{4}\log \kappa_g,
\end{equation}
we get
$$
P_gv+2Q_g=2 Ke^{4v}
$$
Thus, since $0<\kappa_g<8\pi^2$, then the compactness result of Malchiodi\cite{mal} and Druet-Robert\cite{dr} imply \;$\forall m\in \N$, there exists \;$\tilde C_m>0$\; such that
$$
||v||_{C^m(M)}\leq \tilde C_m.
$$
Hence, \;$u\in H^2_Q(M)$\; and \; \eqref{defiv} give the existence of  \;$C_m>0$\; such that
$$
||u||_{C^m(M)}\leq C_m,
$$
thereby ending the proof.
\end{pfn}
\vspace{6pt}

\noindent
\begin{pfn}{ of Theorem \ref{compcri}}\\
The proof is a small modification of the one of Theorem \ref{compsub}. For the sake of completeness, we repeat all the steps. Let $u\in C^{\infty}(M)\cap H^2_Q(M)$ be a minimizer of \;$I$\; on \;$H^2_Q(M)$. Then  as in the proof of Theorem \ref{compsub}, \;$u$ is a minimizer of \;$J$\; on \;$H^2(M)$. Hence \;$u$\; satisfies
$$
P_gu+2Q_g=2\kappa_g \frac{Ke^{4u}}{\int_MKe^{4u}}.
$$
Then, setting 
$$
v=u-\frac{1}{4}\log \int_{M}Ke^{4u}+\frac{1}{4}\log \kappa_g,
$$
we get
$$
P_gv+2Q_g=2 Ke^{4v}.
$$
Thus since \;$\mathcal{F}_{\infty}^0=Crit(\mathcal{F}_K)$, then our compactness theorem in \cite{nd6} imply that \;$\forall m\in \N$\; there exists \;$\tilde C_m>0$ such that
$$
||v||_{C^m(M)}\leq \tilde C_m.
$$
Hence recalling that \;$u\in H^2_Q(M)$, we have there exists \;$C_m>0$\; such that
$$
||u||_{C^m(M)}\leq C_m.
$$
\end{pfn}
\section{Obstacle problem and Moser-Trudinger type inequality}\label{casesphere}
In this section, we discuss some Moser-Trudinger type inequalities related to the Paneitz obstacle problem. In particular, we specialize to the case of the \;$4$-dimensional standard sphere \;$(\mathbb{S}^4, g_{\mathbb{S}^4})$.
\vspace{6pt}

\noindent
We have the following obstacle Moser-Trudinger type inequality.
\begin{pro}
Assuming that  \;$P_g \geq 0$, $\ker P_g=\R$,  then there exists \;$C=C(M, g)>0$\; such that 
$$
\log \int_Me^{4T_g(u)}dV_g\leq C+\frac{1}{8\pi^2}\left<u, u\right>_g, \;\;\;\forall u\in H^2_Q(M).
$$
\end{pro}
\begin{pf}
Clearly \;$u\leq T_g(u)$\; gives
\begin{equation}\label{imp0}
\log \int_Me^{4u}dV_g\leq \log \int_Me^{4T_g(u)}dV_g.
\end{equation}
 Since \;$P_g \geq 0$\; and \;$\ker P_g=\R$, then \;the classical Moser-Trudinger inequality in Proposition \ref{cmt} implies the existence of \; $C=C(M, g)>0$\; such that
\begin{equation}\label{imp1}
\log \int_Me^{4T_g(u)}dV_g\leq C+\frac{1}{8\pi^2}\left<T_g(u), T_g(u)\right>_g.
\end{equation}
Using the definition of \;$T_g$, we get
\begin{equation}\label{imp2}
\left<T_g(u), T_g(u)\right>_g\leq \left<u, u\right>_g.
\end{equation}
Hence combining \eqref{imp0}-\eqref{imp2}, we get
$$
\log \int_Me^{4T_g(u)}dV_g\leq C+\frac{1}{8\pi^2}\left<u, u\right>_g.
$$
\end{pf}
\vspace{8pt}

\noindent
When \;$(M, g)=(\mathbb{S}^4, g_{\mathbb{S}^4})$\; and  \;$K=1$, we have the following sharp obstacle Moser-Trudinger type inequality. 
\begin{thm}\label{sharpomt}
Assuming that \;$(M,  g)=(\mathbb{S}^4, g_{\mathbb{S}^4})$\; and \;$K=1$, then 
$$I\geq 0 \;\;\;\text{on}\;\;\; H^2_Q(M),$$ i.e
\begin{equation}\label{opmt}
\log \int_Me^{4T_g(u)}dV_g\leq \frac{1}{8\pi^2}\left<P_gu, u\right>, \;\;\;\forall u\in H^2_Q(M).
\end{equation}
Moreover equality in \eqref{opmt}  holds if and only if \;$$v:=u-\frac{1}{4}\log \int _Me^{4u}+\frac{1}{4}\log \frac{\kappa_g}{3}$$ is a standard bubble, see \eqref{standard} for its definition.
\end{thm}
\begin{pf}
Since \;$(M, g)=(\mathbb{S}^4, g_{\mathbb{S}^4})$\; and \;$K=1$, then by the classical Moser-Trudinger-Onofiri ineqality in Proposition \ref{cmto}, we have
\begin{equation}\label{imp3}
J\geq 0 \;\;\;\text{on}\;\;\; \;H^2(M)
\end{equation}
and \;\begin{equation}\label{impes}
J(u)=0\;\; \;\text{is equivalent to }\;\;\;v:=u-\frac{1}{4}\log \int _Me^{4u}+\frac{1}{4}\log \frac{\kappa_g}{3}\;\; \;\text{is a standard bubble}.
\end{equation}
 Using Lemma \ref{decreasingformulaop}, we get
\begin{equation}\label{imp4}
I\geq I\circ T_g\; \;\;\text{on}\;\;\; H^2_Q(M).
\end{equation}
Thus, using Lemma \ref{minimaltunnel} and \eqref{imp4}, we have
\begin{equation}\label{imp5}
I\geq J\circ T_g \;\;\;\text{on}\;\;\; H^2_Q(M).
\end{equation}
So, combining \eqref{imp3} and \eqref{imp5}, we get 
\begin{equation}\label{ipositive}
I\geq 0 \;\;\;\text{on}\;\;\; H^2_Q(M).
\end{equation}
Hence, recalling the definition of \;$I$ (see \eqref{eq:defi12}) and \eqref{kgsphere}, we  have \eqref{ipositive} is equivalent to
$$
\log \int_Me^{4T_g(u)}dV_g\leq \frac{1}{8\pi^2}\left<u, u\right>_g, \;\;\;\forall u\in H^2_Q(M).
$$
Suppose \;$$v:=u-\frac{1}{4}\log \int _Me^{4u}+\frac{1}{4}\log \frac{\kappa_g}{3}$$is a standard bubble with\;$u\in H^2_Q(M)$. Then \eqref{impes} implies
\begin{equation}\label{impes1}
J(u)=0
\end{equation}
Thus  \eqref{impes1}, Lemma \ref{minimaltunnel} and the first part (namely \eqref{ipositive})  imply
$$
I(u)=0.
$$
Hence we have the equality case in \eqref{opmt} . Suppose we have the equality case in \eqref{opmt} with \;$u\in H^2_Q(M)$. Then 
\begin{equation}\label{izero}
I(u)=0.
\end{equation}
Thus, using \eqref{ipositive} and \eqref{izero}  we get
\begin{equation}\label{minimum0}
I(u)=\min_{v\in H^2_Q(M)}I(v).
\end{equation}
Using \eqref{minimum0} and Corollary \ref{rigidityminimizeri}, we obtain
\begin{equation}\label{fixedpt}
u=T_g(u).
\end{equation}
So Lemma \ref{minimaltunnel} , \eqref{izero} and \eqref{fixedpt} imply
\begin{equation}\label{jzero}
J(u)=0.
\end{equation}
Hence using \eqref{impes} and \eqref{jzero}, we have \;$v:=u-\frac{1}{4}\log \int _Me^{4u}+\frac{1}{4}\log \frac{\kappa_g}{3}$\; is a standard bubble.
\end{pf}
\vspace{8pt}

\noindent
Theorem \ref{sharpomt} implies the following corollary stating that $Q$-normalized standard bubbles (see \eqref{qnormalized} for their definitions) are fixed points of the obstacle solution map \;$T_g$.
\begin{cor}\label{sfixedt}
Assuming that \;$(M,  g)=(\mathbb{S}^4, g_{\mathbb{S}^4})$\; and \;$w$\; is a \;$Q$-normalized standard bubble (see \eqref{qnormalized} for its definition), then  
$$
T_g(w)=w.
$$
\end{cor}
\begin{pf}
Since \;$w$\; is a \;$Q$-normalized standard bubble, then 
  $$w:=v-\ov{(v)}_Q$$ with \;$v$\; is a standard bubble. Thus, Lemma \ref{minimaltunnel}, Theorem \ref{sharpomt}, and the translation invariant property of \;$J$\; imply
$$
0\leq I(w)\leq J(w)=J(v)=0.
$$
Using again Theorem \ref{sharpomt}, we obtain 
$$
I(w)=\min_{v\in H^2_Q(M)}I(v)
$$
Hence using Corollary \ref{rigidityminimizeri}, we get
$$
w=T_g(w).
$$
\end{pf}

\end{document}